\documentclass[hidelinks,onefignum,onetabnum]{siamart251216}

\usepackage{amsfonts}
\usepackage{graphicx}
\usepackage{epstopdf}
\usepackage{float}
\usepackage{tikz}
\usetikzlibrary{arrows.meta}
\ifpdf
\DeclareGraphicsExtensions{.eps,.pdf,.png,.jpg}
\else
\DeclareGraphicsExtensions{.eps}
\fi

\newsiamremark{remark}{Remark}
\newsiamremark{hypothesis}{Hypothesis}
\crefname{hypothesis}{Hypothesis}{Hypotheses}
\newsiamthm{claim}{Claim}
\newsiamremark{fact}{Fact}
\crefname{fact}{Fact}{Facts}

\headers{The optimal barrier parameter of the power cone}{Rui-Jin Zhang, Wenhao Fu, and Yu-Hong Dai}

\title{The optimal barrier parameter of the three-dimensional power cone\thanks{Submitted to the editors DATE.
		 \funding{ R.J.~Zhang was supported by the National Natural Science Foundation of China (No. 12501432) and the Fundamental Research Funds for the Central Universities (Nos. 050-63253088 and 050-63263071). W.~Fu was supported by the Natural Science Foundation of Jiangsu Province (No. BK20250981). Y.H.~Dai was supported by the National Key  R\&D Program of China (Nos. 2021YFA1000300 and 2021YFA1000301) and the National Natural Science Foundation of China (No. 92473208).}
         }}

\author{Rui-Jin Zhang\thanks{School of Mathematical Sciences and LPMC, Nankai University, Tianjin, 300071, China 
		(\email{zhangrj@nankai.edu.cn}).}
	\and Wenhao Fu\thanks{Corresponding author. School of Mathematical Sciences, Suzhou University of Science and Technology, Suzhou, 215009, China 
		(\email{wenhfu@usts.edu.cn}).}
	\and Yu-Hong Dai\thanks{State Key Laboratory of Mathematical Sciences, Academy of Mathematics and Systems Science, Chinese Academy of Sciences, Beijing 100190, China, and School of Mathematical Sciences,
		University of Chinese Academy of Sciences, Beijing 100049, China 
		(\email{dyh@lsec.cc.ac.cn}).}}
\usepackage{amsopn}
\usepackage{amsmath,amssymb}

\def\K{\mathbb{K}}
\renewcommand{\[}{\begin{equation*}} 
	\renewcommand{\]}{\end{equation*}}

\newcommand{\R}{\mathbb{R}}
\newcommand{\Disc}{\operatorname{Disc}}
\newcommand{\Res}{\operatorname{Res}}
\newcommand{\Kpow}{\K}
\newcommand{\U}{\mathbb{U}}
\newcommand{\W}{\mathbb{W}}
\newcommand{\D}{\mathbb{D}}

\ifpdf
\hypersetup{
  pdftitle={The optimal barrier parameter of the three-dimensional power cone},
  pdfauthor={Rui-Jin Zhang, Wenhao Fu, and Yu-Hong Dai}
}
\fi

\begin{document}

\maketitle

\begin{abstract}
We determine the optimal logarithmically homogeneous self-concordant barrier parameter for the three-dimensional power cone. Although several self-concordant barriers for this cone have been constructed, its optimal barrier parameter had remained unknown. For a power cone with exponent $c\in(0,1)$, we prove that the optimal parameter is $\frac{3}{1+\min\{c,1-c\}}$. Except at $c=1/2$, this value is strictly smaller than the parameter $3-2\min\{c,1-c\}$ of Chares' conjectured barrier. At $c=1/2$, the optimal parameter equals $2$, in agreement with the optimal barrier parameter of the second-order cone. The proof combines a sharp lower bound for arbitrary logarithmically homogeneous self-concordant barriers with an explicit barrier attaining the bound.
\end{abstract}



\section{Introduction}\label{sec1}
A logarithmically homogeneous self-concordant barrier (LHSCB) with parameter $\nu$ leads to short-step interior-point methods with iteration complexity $\mathcal{O}(\sqrt{\nu}\log(1/\varepsilon))$~\cite{NesterovNemirovskii1994,papp2025interior}. The question is then to determine the optimal (smallest) barrier parameter attainable on a given cone. For a symmetric cone associated with a Euclidean Jordan algebra of rank $r$, the log-determinant barrier has parameter $r$, and no self-concordant barrier can have a smaller parameter~\cite{CardosoVieira2006,GulerTuncel1998}. For nonsymmetric cones, however, optimal barrier parameters are known in far fewer cases \cite{fawzi2023optimal,he2026operator}. 

The three-dimensional power cone is a particularly relevant example. Its optimal barrier parameter remained unknown (see~\cite{pirau2025convexity}). The cone is defined by
\[
\Kpow_c
:=
\{(x,y,z)\in\R_+^2\times\R:
|z|\le x^c y^{1-c}\},
\quad 0<c<1.
\]
The case $c=1/2$ is isomorphic to a second-order cone \cite{lin2024generalized}, whereas the remaining members of the family are generally nonsymmetric. 
Several self-concordant barriers have been developed for the power cone. Nesterov~\cite{Nesterov2012} constructed a self-concordant barrier with parameter $4$, thereby placing the power cone within the framework of polynomial-time interior-point methods for nonsymmetric cones.  The construction uses separate logarithmic terms for the two positive coordinates and yields a uniform upper bound on the barrier parameter for all $c\in(0,1)$.

Chares~\cite{Chares2009} reduced the parameter to $3$ by weighting the logarithmic terms according to the exponent.  
At $c=1/2$, however, this bound exceeds the optimal parameter $2$ of the second-order cone~\cite{CardosoVieira2006}. This mismatch led Chares to propose an exponent-dependent candidate with parameter $3-2\min\{c,1-c\}$. Numerical tests~\cite[Section 3.1.2]{Chares2009} suggested that the candidate was self-concordant, but no proof was provided for the entire exponent range. Thus, the value $3-2\min\{c,1-c\}$ remained conjectural rather than an established upper bound.
Chares also proposed a barrier for the $(n,m)$-generalized power cone. Roy and Xiao~\cite{RoyXiao2022} later proved that this barrier is an $(n+1)$-LHSCB. For the standard three-dimensional power cone, corresponding to $n=2$ and $m=1$, their result recovers the parameter $3$ established by Chares.
 
The parameter $3$ can also be obtained from general barrier theory for regular convex cones. Hildebrand~\cite{hildebrand2014canonical} proved that the canonical barrier on an $n$-dimensional regular convex cone is an LHSCB with parameter at most $n$. Since $\Kpow_c$ is three-dimensional, this result again gives the same upper bound $3$.  Hildebrand~\cite{Hildebrand2013} also derived a lower bound depending on the exponent that must be satisfied by every LHSCB on the power cone. The bound varies continuously with the exponent and recovers the optimal barrier parameter $2$ for the second-order cone at $c=1/2$. For $c\ne 1/2$, however, it remains below both Chares' proved parameter $3$ and his conjectured parameter $3-2\min\{c,1-c\}$, leaving a nonzero gap between the known upper and lower bounds.

In this paper, we determine the optimal parameter among all LHSCBs on $\Kpow_c$ by establishing matching lower and upper bounds. For $0<c\leq 1/2$, the lower bound applies to every LHSCB, while the upper bound is attained by an explicit  barrier. The symmetry under $c\leftrightarrow 1-c$ then extends the result to the full range $0<c<1$. The following subsections state the main result and outline its proof.

\subsection{Main result}\label{subsec1.1}
We define
$
\nu_{\mathrm{opt}}(\K)
:=
\inf\{\nu\geq 1:\K\text{ admits a }\nu\text{-LHSCB}\}.
$
By Lemma~\ref{lem:exponent-symmetry}, it suffices to consider $0<c\leq1/2$.
\begin{theorem}[Optimal barrier parameter]\label{thm:main}
Let $0<c\le1/2$. Then the optimal logarithmically homogeneous self-concordant barrier parameter (or optimal barrier parameter for short) of $\Kpow_c$ is
\begin{equation}\label{eq:nu}
\nu_{\rm opt}(\Kpow_c)
=\nu_c
:=
\frac{3}{1+c}.
\end{equation}
Moreover, the optimum is attained by
\begin{equation}\label{eq:Fdef}
\begin{aligned}
F_c(x,y,z)
={}&
-\log x
-\frac{2-c}{1+c}\log y 
-\frac{2-c}{2(c^2-c+1)}
\log\!\left(1-\frac{z^2}{x^{2c}y^{2(1-c)}}\right).
\end{aligned}
\end{equation}
\end{theorem}

\begin{figure}[htbp]
\centering
\includegraphics[width=\textwidth]{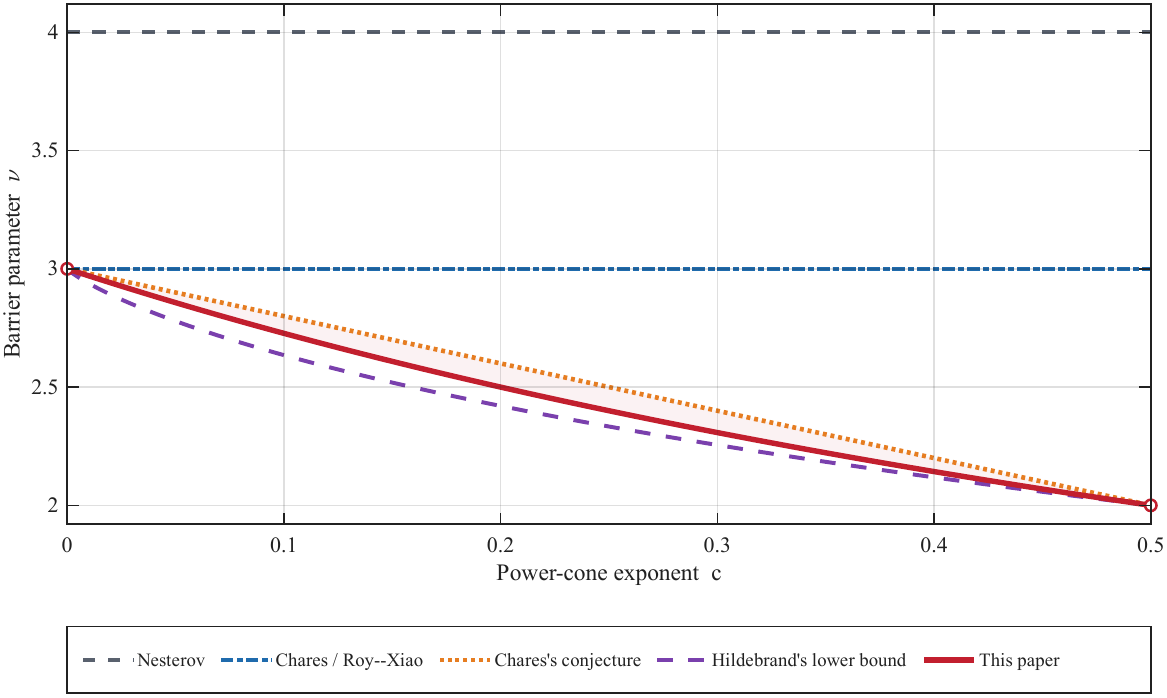}
\caption{Comparison of known, conjectured, and optimal barrier parameters for $0<c\leq1/2$.}
\label{fig:barrier-comparison}
\end{figure}

Figure~\ref{fig:barrier-comparison} shows that our result closes the gap between earlier barrier constructions and lower bounds. In particular, for every $0<c<1/2$, we have
$
\frac{3}{1+c}<3-2c<3.
$
Hence the optimal parameter obtained here is smaller than both the established upper bound $3$  and the parameter $3-2c$ of Chares' conjectured barrier.

\subsection{Proof strategy and organization}\label{subsec:proof-roadmap}
Section~\ref{sec2} introduces the preliminaries needed for the proof. The proof of Theorem~\ref{thm:main} proceeds in two parts. Section~\ref{sec3} establishes the lower bound for every LHSCB by reducing the problem to suitable two-dimensional sections and combining estimates derived from the Dikin ellipsoid property with the projective self-concordance inequality. Section~\ref{sec4} establishes the matching upper bound by verifying that the explicit function $F_c$ is an LHSCB. On an affine section, the self-concordance inequality reduces to the nonnegativity of a binary sextic, which is proved by exact algebraic arguments. Section~\ref{sec5} combines the two bounds to prove the main theorem and extends the result to all $0<c<1$ using Lemma~\ref{lem:exponent-symmetry}. Section~\ref{sec6} concludes the paper. The polynomial factorizations and auxiliary algebraic calculations are given in the appendices. 

\section{Preliminaries}\label{sec2}
This section recalls the basic properties of logarithmically homogeneous self-concordant barriers and presents the algebraic tools used in the upper-bound proof.

\subsection{Logarithmically homogeneous self-concordant barriers}\label{subsec:lhscb-prelim}

Let $\K\subseteq\R^n$ be a regular cone, that is, a closed, convex, and pointed cone with nonempty interior, and let $\K^*$ denote its dual cone. A convex function $F:\operatorname{int}\K\to\R$ is called a $\nu$-logarithmically homogeneous self-concordant barrier ($\nu$-LHSCB), where $\nu\geq1$, if it is three times continuously differentiable with positive-definite Hessian and satisfies
\begin{equation}\label{ineq:selfconcordance}
|D^3F(x)[h,h,h]|
\leq
2\bigl(D^2F(x)[h,h]\bigr)^{3/2}
\end{equation}
for every $x\in\operatorname{int}\K$ and $h\in\R^n$.
In addition, $F$ tends to $+\infty$ whenever an interior sequence
approaches a boundary point of $\K$, and
\[
F(\tau x)=F(x)-\nu\log\tau
\]
for every $x\in\operatorname{int}\K$ and $\tau>0$.
The constant $\nu$ is called the barrier parameter of $F$.

The property of being a $\nu$-LHSCB is preserved under invertible linear transformations \cite[Proposition~2.3.3(i)]{NesterovNemirovskii1994}. In particular, the invertible linear transformation $(x,y,z)\mapsto(y,x,z)$ maps the power cone with exponent $c$ onto the power cone with exponent $1-c$.  The resulting symmetry of the optimal LHSCB parameter is stated in the following lemma, which allows us to restrict the analysis to $0<c\leq 1/2$.
\begin{lemma}\label{lem:exponent-symmetry}
For $0<c<1$, let $T(x,y,z)=(y,x,z)$. Then
\[
T(\Kpow_c)=\Kpow_{1-c}
\quad\text{and}\quad
\nu_{\rm opt}(\Kpow_c)=\nu_{\rm opt}(\Kpow_{1-c}).
\]
More precisely, if $F$ is a $\nu$-LHSCB on $\Kpow_c$, then $F\circ T^{-1}$ is a $\nu$-LHSCB on $\Kpow_{1-c}$.
\end{lemma}

\begin{proof}
For $(x,y,z)\in\Kpow_c$, the inequality
$
|z|\leq x^cy^{1-c}=y^{1-c}x^c
$
shows that $T(x,y,z)=(y,x,z)\in\Kpow_{1-c}$. Since $T^{-1}=T$, the
reverse inclusion follows in the same way. Hence
$
T(\Kpow_c)=\Kpow_{1-c}.
$

Let $F$ be a $\nu$-LHSCB on $\Kpow_c$ and set
$\widetilde F:=F\circ T^{-1}$.
Since $T^{-1}$ is an invertible linear transformation,
it follows from \cite[Proposition~2.3.3(i)]{NesterovNemirovskii1994} that
$\widetilde F$ is a $\nu$-LHSCB on $\Kpow_{1-c}$.
Consequently,
$
\nu_{\rm opt}(\Kpow_{1-c})
\leq
\nu_{\rm opt}(\Kpow_c).
$
Applying the same argument to $T^{-1}=T$ gives the
reverse inequality. Therefore,
$
\nu_{\rm opt}(\Kpow_c)
=
\nu_{\rm opt}(\Kpow_{1-c}).
$
\end{proof}

We also use several standard properties of LHSCBs in the lower-bound argument; see \cite[Chapter~2]{NesterovNemirovskii1994}.
With $H_x:=\nabla^2F(x)$ and $\|h\|_x^2:=h^\top H_xh$, the unit Dikin
ellipsoid satisfies
\begin{equation}\label{eq:dikin-property-background}
\{x+h:\|h\|_x<1\}\subset\operatorname{int}\K.
\end{equation}
Differentiating logarithmic homogeneity gives
\begin{equation}\label{eq:LH-identities-lower}
H_xx=-\nabla F(x),
\quad
\langle-\nabla F(x),x\rangle=\nu.
\end{equation}
Moreover, $-\nabla F(x)\in\operatorname{int}\K^*$.
The conjugate barrier is defined by
\[
F^*(s):=\sup_{x\in\operatorname{int}\K}\{-\langle s,x\rangle-F(x)\}.
\]

The conjugate $F^*$ is a $\nu$-LHSCB on $\K^*$. For $s\in\operatorname{int}\K^*$, let $ H_s^*:=\nabla^2F^*(s)$ and  $\|u\|_{s,*}^2:=u^\top H_s^*u$. Its unit Dikin ellipsoid satisfies 
\begin{equation}\label{eq:dual-dikin-property-background} \{s+u:\|u\|_{s,*}<1\} \subset\operatorname{int}\K^*. \end{equation}
At $s=-\nabla F(x)$, the conjugacy relations give
\begin{equation}\label{eq:conjugate-hessian-lower}
\nabla F^*(s)=-x,
\quad
\nabla^2F^*(s)=H_x^{-1}.
\end{equation}

Hildebrand~\cite[Theorem~2]{Hildebrand2022} established an equivalence between logarithmically homogeneous self-concordant functions on cones and projectively self-concordant functions on their affine sections.

\begin{lemma}\label{lem:projective}
Let $\D\subset\R^n$ be open and convex, let $\nu\ge2$, and let $f$ be three times continuously differentiable on $\D$. Consider the $\nu$-logarithmically homogeneous function associated with $f$, given by
\[
\widehat F(t,w)=\nu\bigl(-\log t+f(t^{-1}w)\bigr),
\quad
(t,w)\in\{(t,w):t>0,\ t^{-1}w\in \D\}.
\]
Define the bilinear form and the cubic form
\begin{align*}
G&:=\nabla^2f-\nabla f\nabla f^\top,\\
\mathcal T_f[h,h,h]
&:=D^3f[h,h,h]-6D^2f[h,h]Df[h]+4(Df[h])^3.
\end{align*}
Then $\widehat F$ has positive-definite Hessian and satisfies the standard self-concordance inequality \eqref{ineq:selfconcordance} if and only if $G\succ0$ and
\begin{equation}\label{eq:projectiveSC}
|\mathcal T_f[h,h,h]|
\le
2\gamma\,G[h,h]^{3/2}
\quad
\text{for all }h,
\end{equation}
where $\gamma=(\nu-2)/\sqrt{\nu-1}$.
\end{lemma}

\subsection{Algebraic tools}\label{subsec:poly-prelim}
We collect the standard algebraic tools used in the upper-bound proof; see \cite{BasuPollackRoy2006,FaroukiRajan1988,HarrisMartin1987}. For polynomials $p$ and $q$ of degrees $m$ and $n$, $\Res_\xi(p,q)$ denotes their resultant with respect to $\xi$. It is the determinant of the $(m+n)\times(m+n)$ Sylvester matrix whose rows are the coefficient vectors, ordered by descending powers of $\xi$, of $\xi^{n-1}p,\ldots,p$ and $\xi^{m-1}q,\ldots,q$; see \cite[Chapter~4]{BasuPollackRoy2006}. For $p(\xi)=\sum_{j=0}^d a_j\xi^j$, write $[\xi^j]p:=a_j$. If $p$ has degree $d$ and leading coefficient $a_d\ne0$, its discriminant is
\begin{equation}\label{eq:disc-definition}
\Disc_\xi p:=(-1)^{d(d-1)/2}\frac{\Res_\xi(p,p')}{a_d}.
\end{equation}

\begin{lemma}\label{lem:disc-background}
A polynomial has zero discriminant if and only if it has a multiple complex root. If $p$ is a real cubic with nonzero discriminant, then
\[
\Disc p<0
\quad\Longleftrightarrow\quad
p\text{ has exactly one real root}.
\]
\end{lemma}

\begin{lemma}\label{lem:root-continuation}
Let $\{p_\theta\}$ be a continuous family of real polynomials of fixed degree, where $\theta$ ranges over a connected parameter set. If $\Disc p_\theta\ne0$ for every $\theta$, then the number of real roots of $p_\theta$ is independent of $\theta$.
\end{lemma}
\begin{proof}
Polynomial roots depend continuously on the coefficients \cite{HarrisMartin1987}. All roots remain simple since $\Disc p_\theta\ne0$. Hence, under a sufficiently small change of $\theta$, each real root remains real and each nonreal conjugate pair remains nonreal. The number of real roots is therefore locally constant. Since the parameter set is connected, it is constant throughout the set. 
\end{proof}

We use the following form of Sturm's theorem to count real roots on $(-\infty,+\infty)$.

\begin{lemma}\label{lem:sturm-background}
 Let $p$ be a nonconstant real polynomial. For polynomials $f$ and nonzero $g$, let $\operatorname{rem}(f,g)$ denote the remainder obtained by Euclidean division of $f$ by $g$. The Sturm chain of $p$ is defined by $S_0=p$, $S_1=p'$, and $S_{j+1}=-\operatorname{rem}(S_{j-1},S_j)$ for $j\ge1$ until a zero remainder is obtained. Then the number of distinct real roots of $p$ is $V(-\infty)-V(+\infty)$, where $V$ counts sign changes in the Sturm chain, with the signs at $\pm\infty$ determined by the leading coefficients and degrees of its members. Multiplying any member of the chain by a positive constant does not change the root count.
\end{lemma}

For $m\ge0$, the Bernstein basis of degree $m$ on $[0,1]$ is
\[
B_i^{(m)}(v):=\binom{m}{i}v^i(1-v)^{m-i},
\quad 0\le i\le m.
\]
These basis functions are nonnegative on $[0,1]$ and satisfy $ \sum_{i=0}^m B_i^{(m)}(v)=1. $ Hence, if a polynomial $p$ of degree at most $m$ is represented in the Bernstein basis as
$
p(v)=\sum_{i=0}^m b_i B_i^{(m)}(v),
$
then $b_i\ge0$ for all $i$ implies $p(v)\ge0$ on $[0,1]$, while $b_i>0$ for all $i$ implies $p(v)>0$ on $[0,1]$. The following identity converts the power basis $v^k$ into the Bernstein basis of degree $m$:
\begin{equation}\label{eq:power-to-bernstein}
v^k=\sum_{i=k}^m
\frac{\binom{i}{k}}{\binom{m}{k}}B_i^{(m)}(v),
\quad 0\le k\le m.
\end{equation}

\section{The lower bound}\label{sec3}
This section establishes the lower bound for every LHSCB on $\Kpow_c$ needed for Theorem~\ref{thm:main}.
\subsection{Auxiliary results}\label{subsec3.1}
The proof relies on four auxiliary results, which are presented below in the order in which they are used. The first result is the  partial-minimization property of self-concordant barriers.
\begin{lemma}\label{lem:partialmin} 
Let $\K\subset\U\times\W$ be a regular cone, and let $\bar\K:=\Pi_{\U}(\K)$ be its projection onto $\U$. Assume that $\bar\K$ is also a regular cone. For $u\in\operatorname{int}\bar\K$, define 
\[
\K_u:=\{w\in\W:(u,w)\in\operatorname{int}\K\}. 
\]
Let $\Phi$ be a $\nu$-LHSCB on $\K$. Assume that $w\mapsto\Phi(u,w)$ attains its minimum over $\K_u$ at $\zeta(u)$, and define 
\[
\bar\Phi(u):=\Phi(u,\zeta(u)).
\]
Suppose further that whenever $u_k\to\bar u\in\partial\bar\K$, the points $(u_k,\zeta(u_k))$ converge to a point of $\partial\K$. Then $\bar\Phi$ is a $\nu$-LHSCB on $\bar\K$. 
\end{lemma} 
\begin{proof} 
Since $\nabla^2\Phi$ is positive definite, the minimizer $\zeta(u)$ is unique and satisfies 
\begin{equation}\label{eq:partialmin-stationarity}
D_w\Phi(u,\zeta(u))=0.
\end{equation}
Since the $ww$ block of $\nabla^2\Phi$ is also positive definite, the implicit-function theorem implies that $\zeta$ is twice continuously differentiable. Differentiating \eqref{eq:partialmin-stationarity} gives
\begin{equation}\label{eq:partialmin-hessian}
\nabla^2\bar\Phi(u)
=
\nabla^2_{uu}\Phi(u,\zeta(u))-
\nabla^2_{uw}\Phi(u,\zeta(u))
\bigl[\nabla^2_{ww}\Phi(u,\zeta(u))\bigr]^{-1}
\nabla^2_{wu}\Phi(u,\zeta(u)).
\end{equation}
Since $\Phi$ is three times continuously differentiable and $\zeta$ is twice continuously differentiable, the right-hand side of \eqref{eq:partialmin-hessian} is continuously differentiable. Hence $\bar\Phi$ is three times continuously differentiable.

Fix $u$ and $h\in\U$, and set $X(t):=(u+th,\zeta(u+th))$. Differentiating \eqref{eq:partialmin-stationarity} along $X(t)$ gives
\[
D^2\Phi(X(t))[X'(t),(0,\omega)]=0
\quad\text{for every }\omega\in\W.
\]
Since $X''(t)\in\{0\}\times\W$, it follows that $D^2\Phi(X(t))[X'(t),X''(t)]=0$. The stationarity condition also gives $D\Phi(X(t))[X''(t)]=0$. Therefore, differentiating $\bar\Phi(u+th)=\Phi(X(t))$ and evaluating at $t=0$ yields
\[
D^2\bar\Phi(u)[h,h]=D^2\Phi(X(0))[v,v],\quad
D^3\bar\Phi(u)[h,h,h]=D^3\Phi(X(0))[v,v,v],
\]
where $v:=X'(0)$. Since $h\ne0$ implies $v\ne0$, the first identity gives $\nabla^2\bar\Phi\succ0$, while the self-concordance of $\Phi$ gives
\[
|D^3\bar\Phi(u)[h,h,h]|
\le 2\bigl(D^2\bar\Phi(u)[h,h]\bigr)^{3/2}.
\]

For $t>0$, the conic structure gives $\K_{tu}=t\K_u$. Therefore \[ \bar\Phi(tu) =\min_{w\in\K_u}\Phi(tu,tw) =\bar\Phi(u)-\nu\log t, \] so $\bar\Phi$ is $\nu$-logarithmically homogeneous. Finally, if $u_k\to\bar u\in\partial\bar\K$, the boundary assumption and the barrier property of $\Phi$ give $\bar\Phi(u_k)=\Phi(u_k,\zeta(u_k))\to+\infty$. Thus $\bar\Phi$ is a $\nu$-LHSCB on $\bar\K$.
\end{proof}

The next estimate uses the Dikin ellipsoid properties to obtain two-sided bounds at both endpoints of the chord.
\begin{lemma}\label{lem:dikinchord} 
Let $\K$ be a regular two-dimensional cone, let $\Phi$ be a $\nu$-LHSCB on $\K$, and let $x\in\operatorname{int}\K$. Choose $h$ such that $x$ and $h$ are linearly independent, and suppose that
\[
y_-:=x-\delta_-h\in\partial\K\setminus\{0\},
\quad
y_+:=x+\delta_+h\in\partial\K\setminus\{0\},
\quad
\delta_\pm>0.
\]
Introduce the normalized derivatives
\begin{equation}\label{eq:chord-abq}
a:=\frac{1}{\nu}D\Phi(x)[h], \quad b:=\frac{1}{\nu}D^2\Phi(x)[h,h], \quad q:=b-a^2.
\end{equation} 
Then $q>0$ and $\nu\ge2$. Moreover, with $\varrho:=\nu-1$,
\begin{align}
\frac{1}{\sqrt{\varrho}} \le \frac{a+\delta_-^{-1}}{\sqrt{q}} \le \sqrt{\varrho},\label{eq:chord-left} \\
\frac{1}{\sqrt{\varrho}} \le \frac{-a+\delta_+^{-1}}{\sqrt{q}} \le \sqrt{\varrho}. \label{eq:chord-right}
\end{align} 
\end{lemma} 
\begin{proof} 
Let $H:=\nabla^2\Phi(x)$. Logarithmic homogeneity gives $Hx=-\nabla\Phi(x)$ and $x^{\top}Hx=\nu$. Hence, in the basis $\{x,h\}$, \begin{equation}\label{eq:Hchordmatrix} 
\frac{1}{\nu}H = \begin{pmatrix} 1 & -a\\ -a & b \end{pmatrix}. 
\end{equation} 
Since $H\succ0$, its determinant is positive, and therefore $q=b-a^2>0$. We first prove \eqref{eq:chord-right}. Set $A_+:=1-a\delta_+$. Since $s:=-\nabla\Phi(x)\in\operatorname{int}\K^*$ and $y_+\in\partial\K\setminus\{0\}$, we have
\[ 
\langle s,y_+\rangle = \nu A_+ >0.
\]
Thus $A_+>0$. The unit Dikin ellipsoid at $x$ is contained in $\operatorname{int}\K$ by \eqref{eq:dikin-property-background}. Since $\lambda y_+\in\partial\K$ for every $\lambda>0$, it follows that
\begin{equation}\label{eq:primal-dikin-bound}
\|\lambda y_+-x\|_x\ge1 \quad \text{for every }\lambda>0. 
\end{equation}
In the basis $\{x,h\}$, 
\[
\lambda y_+-x = (\lambda-1)x+\lambda\delta_+h.
\] 
Using \eqref{eq:Hchordmatrix} and $b=q+a^2$, we obtain \begin{align*} 
\frac{1}{\nu}\|\lambda y_+-x\|_x^2 &= (\lambda-1)^2 -2a\delta_+\lambda(\lambda-1) +b\delta_+^2\lambda^2\\ &= 1-2A_+\lambda +\bigl(A_+^2+q\delta_+^2\bigr)\lambda^2. 
\end{align*}
Since $A_+>0$, the minimizer over $\lambda>0$ is $ \lambda_* = \frac{A_+}{A_+^2+q\delta_+^2}>0$. Substitution gives 
\begin{equation}\label{eq:primal-dikin-min} \min_{\lambda>0} \|\lambda y_+-x\|_x^2 = \nu\, \frac{q\delta_+^2}{A_+^2+q\delta_+^2} \ge1, \end{equation} 
where the inequality follows from the preceding bound \eqref{eq:primal-dikin-bound}. Therefore,
\begin{equation}\label{eq:primal-chord-ineq} A_+^2 \le (\nu-1)q\delta_+^2. 
\end{equation} 
The reverse estimate is obtained from the Dikin ellipsoid of the
conjugate barrier $\Phi^*$ on $\K^*$. Since $y_+$ is a nonzero boundary point of $\K$, there exists $\eta_+\in\partial\K^*\setminus\{0\}$ such that $\langle\eta_+,y_+\rangle=0$. Since $x\in\operatorname{int}\K$, $\langle\eta_+,x\rangle>0$, so we may normalize $\eta_+$ such that $\langle\eta_+,x\rangle=1$. It follows from $y_+=x+\delta_+h$ that $\langle\eta_+,h\rangle=-\delta_+^{-1}$. At $s=-\nabla\Phi(x)$, the conjugacy relation gives $\nabla^2\Phi^*(s)=H^{-1}$. In the dual basis associated with $\{x,h\}$,
\[
s=(\nu,-\nu a),
\quad
\eta_+=\left(1,-\delta_+^{-1}\right).
\]
Using these coordinate representations gives
\begin{align*} 
\|\mu\eta_+-s\|_{s,*}^2 = \frac{1}{\nu q\delta_+^2} \Bigl[ \bigl(A_+^2+q\delta_+^2\bigr)\mu^2 -2\nu q\delta_+^2\mu +\nu^2q\delta_+^2 \Bigr]. 
\end{align*} 
Since $\mu\eta_+\in\partial\K^*$ for every $\mu>0$, the Dikin ellipsoid property of $\Phi^*$ gives
\[
\|\mu\eta_+-s\|_{s,*}\ge1 \quad \text{for every }\mu>0. 
\] 
The preceding quadratic is minimized at $\mu_* = \frac{\nu q\delta_+^2}{A_+^2+q\delta_+^2}>0$. Consequently, 
\begin{equation}\label{eq:dual-dikin-min} 
\min_{\mu>0} \|\mu\eta_+-s\|_{s,*}^2 = \nu\, \frac{A_+^2}{A_+^2+q\delta_+^2} \ge1. 
\end{equation}
Hence \begin{equation}\label{eq:dual-chord-ineq} q\delta_+^2 \le (\nu-1)A_+^2. \end{equation} Since $A_+>0$, $q>0$, and $\delta_+>0$, \eqref{eq:primal-chord-ineq} implies $\nu-1>0$. Combining \eqref{eq:primal-chord-ineq} and \eqref{eq:dual-chord-ineq} gives $1\le(\nu-1)^2$. Therefore, $\nu\ge2$. Moreover, the same two inequalities give 
\[ 
\frac{1}{\sqrt{\nu-1}} \le \frac{A_+}{\delta_+\sqrt{q}} \le \sqrt{\nu-1}.
\] 
Since $\frac{A_+}{\delta_+} = -a+\delta_+^{-1},$ this proves \eqref{eq:chord-right}. Finally, replacing $h$ by $-h$ interchanges the two endpoints, replaces $a$ by $-a$, and leaves $q$ unchanged. Applying \eqref{eq:chord-right} to $-h$ therefore yields \eqref{eq:chord-left}.
\end{proof}

The following lemma expresses the projective self-concordance inequality from Lemma~\ref{lem:projective} in terms of the quantities introduced in Lemma~\ref{lem:dikinchord}.
\begin{lemma}\label{lem:radialcubic}
Under the assumptions and notation of Lemma~\ref{lem:dikinchord}, define
$
d:=\frac{1}{\nu}D^3\Phi(x)[h,h,h].
$
Then
\begin{equation}\label{eq:radialcubic}
\left|d-6ab+4a^3\right|
\leq
2\frac{\nu-2}{\sqrt{\nu-1}}\,q^{3/2}.
\end{equation}
\end{lemma}

\begin{proof}
Let $\D:=(-\delta_-,\delta_+)$ and define $f(\tau):=\nu^{-1}\Phi(x+\tau h)$ for $\tau\in \D$. Since $\K$ is a regular two-dimensional cone, its nonzero boundary points $y_-$ and $y_+$ generate $\K$. In the basis $\{x,h\}$, the interior of $\K$ is $\{t x+w h: t>0, \, t^{-1}w \in \D\}$. By logarithmic homogeneity,
\[
\Phi(tx+wh)
=
\nu\bigl(-\log t+f(t^{-1}w)\bigr).
\]
Thus, in these coordinates, $\Phi$ has the form required in Lemma~\ref{lem:projective}, which applies since $\nu\ge2$ by Lemma~\ref{lem:dikinchord}. As $f'(0)=a$, $f''(0)=b$, and $f'''(0)=d$, we obtain $G[1,1]=b-a^2=q$ and $\mathcal T_f[1,1,1]=d-6ab+4a^3$ at $\tau=0$. Inequality~\eqref{eq:projectiveSC} then gives \eqref{eq:radialcubic}.
\end{proof}

We also need an estimate for a mixed component of a symmetric trilinear form.

\begin{lemma}\label{lem:mixedcubic}
Let $G$ be an inner product on a two-dimensional real vector space, and let $\mathcal T$ be a symmetric trilinear form satisfying
\[
|\mathcal T[v,v,v]|
\leq
M\,G[v,v]^{3/2}
\quad
\text{for every }v.
\]
If $E_1$ and $E_2$ are $G$-orthonormal, then
\[
|\mathcal T[E_1,E_2,E_2]|
\leq
M.
\]
\end{lemma}

\begin{proof}
Let $v(\theta):=\cos\theta\,E_1+\sin\theta\,E_2$ and $T(\theta):=\mathcal T[v(\theta),v(\theta),v(\theta)]$. Then $T$ is a trigonometric polynomial of degree at most three, while $G[v(\theta),v(\theta)]=1$ implies $\|T\|_\infty\le M$. Hence the Bernstein inequality for trigonometric polynomials~\cite[Chapter~X, Section~3]{Zygmund2003} gives $\|T'\|_\infty\le 3\|T\|_\infty \leq 3M$. At $\theta=\pi/2$, symmetry and trilinearity yield
\[
T'(\pi/2)=-3\mathcal T[E_1,E_2,E_2],
\]
and the conclusion follows.
\end{proof}

\subsection{Proof of the lower bound}\label{subsec:3.2}

The preceding estimates yield the following lower bound for the barrier parameter.

\begin{theorem}\label{thm:lowerbound}
Let $0<c\leq1/2$. Every $\nu$-LHSCB on $\Kpow_c$ satisfies
\begin{equation}\label{eq:sharp-lower}
\nu\geq\frac{3}{1+c}.
\end{equation}
\end{theorem}

\begin{proof}
Let $\Phi$ be a $\nu$-LHSCB on $\Kpow_c$.
If $c=1/2$, let
$
L:=\operatorname{span}\{(1,1,0),(0,0,1)\}.
$
Since $L\cap\operatorname{int}\Kpow_{1/2}\neq\varnothing$, the
restriction of $\Phi$ to $\Kpow_{1/2}\cap L$ is again a $\nu$-LHSCB~\cite[Proposition~2.3.3(i)]{NesterovNemirovskii1994}.
At $x=(1,1,0)$ in the direction $h=(0,0,1)$, both chord endpoints are
nonzero boundary points. Lemma~\ref{lem:dikinchord} therefore gives
$\nu\geq2=3/(1+c)$ when $c=1/2$. It remains to consider $0<c<1/2$. Suppose, to the contrary, that
\begin{equation}\label{eq:contradict-nu}
\nu<\frac{3}{1+c}.
\end{equation}

For $x,\,y>0$, the function $z\mapsto\Phi(x,y,z)$ is strictly convex on
$(-x^cy^{1-c},x^cy^{1-c})$ and tends to $+\infty$ at both endpoints.
It therefore has a unique interior minimizer. Define
\[
\bar\Phi(x,y):=\min_{|z|<x^cy^{1-c}}\Phi(x,y,z).
\]
If $(x_j,y_j)\to(\bar x,\bar y)\in\partial\R_+^2$ and $z_j$ is the
corresponding minimizing coordinate, then
$|z_j|<x_j^cy_j^{1-c}\to0$. Hence
$(x_j,y_j,z_j)\to(\bar x,\bar y,0)\in\partial\Kpow_c$.
Lemma~\ref{lem:partialmin} therefore shows that $\bar\Phi$ is a
$\nu$-LHSCB on $\R_+^2$.

For every $(\hat x,\hat y,\hat z)\in\operatorname{int}\Kpow_c$, let $s=\hat x+\hat y$, $r=\hat x/s$, and $z=\hat z/s$. Then $s>0$ and
$
(\hat x,\hat y,\hat z)=s(r,1-r,z).
$
By logarithmic homogeneity,
\[
\Phi(\hat x,\hat y,\hat z)=\Phi(r,1-r,z)-\nu\log s.
\]
Thus the restriction of $\Phi$ to the section $x+y=1$ determines $\Phi$ on the whole cone. We therefore set  $f(r,z):=\nu^{-1}\Phi(r,1-r,z)$. Let $\zeta(r)$ be the unique minimizer
of $z\mapsto f(r,z)$, so
\begin{equation}\label{eq:fz0lower}
f_z(r,\zeta(r))=0.
\end{equation}
All derivatives below are evaluated at $(r,\zeta(r))$. Write
\begin{equation}\label{eq:lower-derivatives}
p:=f_r,\quad
a_0:=f_{rr}-p^2,\quad
m:=f_{rz},\quad
k:=f_{zz},
\end{equation}
and define
\begin{equation}\label{eq:lower-G}
G:=\nabla^2f-\nabla f\nabla f^\top
=
\begin{pmatrix}
a_0&m\\
m&k
\end{pmatrix}.
\end{equation}

Let $x_0:=(r,1-r,\zeta(r))$. For a nonzero $v=(v_r,v_z)$, set
$\widehat v:=(v_r,-v_r,v_z)$. Since $x_0$ and $\widehat v$ are linearly
independent, logarithmic homogeneity gives
\[
\frac{1}{\nu}
\begin{pmatrix}
x_0^\top\\
\widehat v^\top
\end{pmatrix}
\nabla^2\Phi(x_0)
\begin{pmatrix}
x_0&\widehat v
\end{pmatrix}
=
\begin{pmatrix}
1&-Df[v]\\
-Df[v]&D^2f[v,v]
\end{pmatrix}
\succ0.
\]
Its Schur complement gives
$G[v,v]=D^2f[v,v]-(Df[v])^2>0$. Thus $G\succ0$, and in particular
$k>0$.

The implicit function theorem applied to \eqref{eq:fz0lower} gives
$\zeta'=-m/k$.  Thus the tangent vector to
$r\mapsto(r,\zeta(r))$ is $e_r:=(1,-m/k)$.
With $e_z:=(0,1)$, define $\widehat q:=G[e_r,e_r]
=a_0-\frac{m^2}{k}>0$.
Then $G[e_r,e_z]=0$. Moreover, the function $g(r):=\nu^{-1}\bar\Phi(r,1-r)=f(r,\zeta(r))$ satisfies
\begin{equation}\label{eq:gmetriclower}
g'=p,\quad g''-(g')^2=\widehat q.
\end{equation}

Apply Lemma~\ref{lem:dikinchord} to $\bar\Phi$ at $(r,1-r)$ in the
direction $(1,-1)$. Set $\varrho:=\nu-1\geq1$ and
\begin{equation}\label{eq:LUdef}
\ell_-:=\frac{p+r^{-1}}{\sqrt{\widehat q}},
\quad
\ell_+:=\frac{-p+(1-r)^{-1}}{\sqrt{\widehat q}}.
\end{equation}
Then
\begin{equation}\label{eq:LUbounds}
\frac{1}{\sqrt{\varrho}}
\leq \ell_-,\,\ell_+
\leq \sqrt{\varrho}.
\end{equation}
The definitions also give
\begin{align}
\ell_-+\ell_+
&=
\frac{1}{r(1-r)\sqrt{\widehat q}},
\label{eq:LU-sum}\\
\frac{p}{\sqrt{\widehat q}}
+\frac{c}{r\sqrt{\widehat q}}
-\frac{1-c}{(1-r)\sqrt{\widehat q}}
&=
c\ell_--(1-c)\ell_+.
\label{eq:LU-combination}
\end{align}

Define the symmetric trilinear form
\begin{align}
\mathcal T_f[u,v,w]
:={}&
D^3f[u,v,w]
-2Df[u]D^2f[v,w]
-2Df[v]D^2f[u,w]\notag\\
&-2Df[w]D^2f[u,v]
+4Df[u]Df[v]Df[w].
\label{eq:polar-cubic-lower}
\end{align}
For any nonzero direction $v$, Lemma~\ref{lem:radialcubic}
applied to the restriction of $\Phi$ to
$\Kpow_c\cap\operatorname{span}\{x_0,\widehat v\}$ gives
\begin{equation}\label{eq:section-cubic-bound}
|\mathcal T_f[v,v,v]|
\leq
2\gamma G[v,v]^{3/2},
\quad
\gamma:=\frac{\varrho-1}{\sqrt{\varrho}}.
\end{equation}
The vectors
$E_1:=e_r/\sqrt{\widehat q}$ and $E_2:=e_z/\sqrt{k}$ are
$G$-orthonormal. Lemma~\ref{lem:mixedcubic} therefore gives
$|\mathcal T_f[E_1,E_2,E_2]|\leq2\gamma$.

Since $k(r)=f_{zz}(r,\zeta(r))$, we have
$k'=D^3f[e_r,e_z,e_z]$. Moreover,
$Df[e_z]=0$, $Df[e_r]=p$, and $D^2f[e_z,e_z]=k$. Hence
$\mathcal T_f[e_r,e_z,e_z]=k'-2pk$. By trilinearity,
\begin{equation}\label{eq:Bboundlower}
\beta:=\mathcal T_f[E_1,E_2,E_2]=
\frac{k'/k-2p}{\sqrt{\widehat q}}
\quad\text{satisfies}\quad
|\beta|
\leq 2\gamma=
\frac{2(\varrho-1)}{\sqrt{\varrho}}.
\end{equation}

Define $w(r):=r^c(1-r)^{1-c}$. We first prove that
$
\kappa(r):=k(r)w(r)^2
$
tends to zero as $r\to1^{-}$. Fix $r_0\in(0,1)$ and define
\[
\sigma(r):=
\int_{r_0}^{r}\sqrt{\widehat q(t)}\,\text{d}t.
\]
Then $\text{d}\sigma/\text{d}r=\sqrt{\widehat q}$, and hence
$
\frac{\text{d}}{\text{d}\sigma}
=
\frac{1}{\sqrt{\widehat q}}
\frac{\text{d}}{\text{d}r}.
$
Since $w'/w=c/r-(1-c)/(1-r)$, equations
\eqref{eq:LU-combination} and \eqref{eq:Bboundlower} yield
\begin{align}
\frac{\text{d}}{\text{d}\sigma}\log\kappa
&=
\beta
+2\left(
\frac{p}{\sqrt{\widehat q}}
+\frac{c}{r\sqrt{\widehat q}}
-\frac{1-c}{(1-r)\sqrt{\widehat q}}
\right)\notag\\
&=
\beta+2c\ell_--2(1-c)\ell_+\notag\\
&\leq
\frac{2}{\sqrt{\varrho}}
\bigl((1+c)\varrho-(2-c)\bigr).
\label{eq:drift-bound}
\end{align}
Assumption \eqref{eq:contradict-nu} is equivalent to
$\varrho<(2-c)/(1+c)$. Thus
\[
\eta:=
\frac{2}{\sqrt{\varrho}}
\bigl((2-c)-(1+c)\varrho\bigr)>0,
\]
and \eqref{eq:drift-bound} gives
\begin{equation}\label{ineq:diff}
\frac{\text{d}}{\text{d}\sigma}\log\kappa\leq-\eta.
\end{equation}

Equation \eqref{eq:LU-sum} and inequality \eqref{eq:LUbounds} imply
\[
\sqrt{\widehat q}
=
\frac{1}{r(1-r)(\ell_-+\ell_+)}
\geq
\frac{1}{2\sqrt{\varrho}\,r(1-r)}.
\]
Therefore
\[
\sigma(r)-\sigma(r_0)
\geq
\frac{1}{2\sqrt{\varrho}}
\int_{r_0}^{r}\frac{\text{d}t}{t(1-t)}
\longrightarrow+\infty
\quad \text{as } r\to1^{-}.
\]
Integrating the differential inequality \eqref{ineq:diff} gives
\begin{equation}\label{eq:Ktozero-lower}
\kappa(r)
\leq
\kappa(r_0)
\exp\!\bigl(-\eta(\sigma(r)-\sigma(r_0))\bigr)
\longrightarrow0
\quad \text{as } r\to1^{-}.
\end{equation}

We now prove that $\kappa$ has a strictly positive lower bound. Fix
$r\in(0,1)$ and consider
\[
\Kpow_c\cap
\operatorname{span}\{(r,1-r,0),(0,0,1)\}.
\]
At the point $x_0=(r,1-r,\zeta(r))$, consider the direction $e_3:=(0,0,1)$. The distances to the boundary in the directions $-e_3$ and $e_3$ are  $\delta_-:=w(r)+\zeta(r)$ and $\delta_+:=w(r)-\zeta(r)$, respectively. By \eqref{eq:fz0lower}, the normalized first
derivative is zero, while the normalized second derivative is $k$.
Lemma~\ref{lem:dikinchord} gives
\[
\frac{1}{\sqrt{\varrho}}
\leq
\delta_\pm\sqrt{k}
\leq
\sqrt{\varrho}.
\]
Multiplying the two lower bounds and using
$\delta_-\delta_+=w(r)^2-\zeta(r)^2\leq w(r)^2$ yields
\begin{equation}\label{eq:Klowerbound}
\kappa(r)=k(r)w(r)^2
\geq
\frac{1}{\varrho}>0
\quad \text{for }
0<r<1.
\end{equation}
This contradicts \eqref{eq:Ktozero-lower}. Thus
\eqref{eq:contradict-nu} is impossible, and
\eqref{eq:sharp-lower} follows.
\end{proof}

\section{The matching upper bound}\label{sec4}
In this section, we show that the function $F_c$ defined in \eqref{eq:Fdef} has barrier parameter $\nu_c$, which yields an upper bound on the optimal barrier parameter of the power cone.

\subsection{Logarithmic homogeneity and the barrier property}\label{subsec:4.1}
We begin by verifying its logarithmic homogeneity and barrier property.
For every $\tau>0$, we have
$
\frac{(\tau z)^2}{(\tau x)^{2c}(\tau y)^{2(1-c)}}
=
\frac{z^2}{x^{2c}y^{2(1-c)}}.
$
Consequently,
\[
F_c(\tau x,\tau y,\tau z)
=
F_c(x,y,z)
-\left(1+\frac{2-c}{1+c}\right)\log\tau
=
F_c(x,y,z)-\nu_c\log\tau.
\]
Hence $F_c$ is $\nu_c$-logarithmically homogeneous.

Let $\{(x_j,y_j,z_j)\}\subset\operatorname{int}\Kpow_c$ be a sequence converging to a boundary point $(\bar x,\bar y,\bar z)\in\partial\Kpow_c$. If $\bar x=0$ or $\bar y=0$, at least one of the first two logarithmic terms in $F_c$ tends to $+\infty$, while the other is bounded below. Moreover, the last term in $F_c$ is nonnegative. If $\bar x,\bar y>0$, the boundary condition yields $|\bar z|=\bar x^c\bar y^{1-c}$. Consequently, the last term tends to $+\infty$, while the first two logarithmic terms remain bounded. Thus, $F_c$ has the barrier property. It remains to establish the positive definiteness of $\nabla^2F_c$ and the self-concordance inequality.

\subsection{Reduction to the affine section \texorpdfstring{$y=1$}{yequal1}}\label{subsec:4.2}
The endpoint $c=1/2$ will be handled directly in the proof of
Theorem~\ref{thm:upperbound}.
Assume from now on that $0<c<1/2$ and introduce
$
e:=1-2c\in(0,1).
$
Then
\begin{equation}\label{eq:nu-e}
\nu_c=\frac{6}{3-e},
\quad
\gamma
=
\frac{\nu_c-2}{\sqrt{\nu_c-1}}
=
\frac{2e}{\sqrt{9-e^2}}.
\end{equation}
Take the affine section $y=1$. Every interior point has $y>0$ and can be written uniquely as
\[
(x,y,z)=y(X,1,Z),
\quad
X:=\frac{x}{y},
\quad
Z:=\frac{z}{y}.
\]
By logarithmic homogeneity,
\begin{equation}\label{eq:Fprojective}
F_c(x,y,z)
=
\nu_c\bigl[-\log y+f_e(X,Z)\bigr],
\end{equation}
where $f_e$ is defined by
\begin{equation}\label{eq:fdef}
f_e(X,Z)
:=
-\frac{3-e}{6}\log X
-\frac{9-e^2}{6(e^2+3)}
\log\!\left(1-\frac{Z^2}{X^{1-e}}\right).
\end{equation}
The domain of $f_e$ is
$
\D_e
:=
\{(X,Z):X>0,\ |Z|<X^{(1-e)/2}\}.
$
For every $\lambda>0$,
\[
f_e(\lambda X,\lambda^{(1-e)/2}Z)
=
f_e(X,Z)-\frac{3-e}{6}\log\lambda.
\]
Since $(1-e)/2=c$, the corresponding linear map is
\[
S_\lambda(X,Z):=(\lambda X,\lambda^cZ).
\]
The additive term on the right is independent of $(X,Z)$. Consequently, for every direction $h$,
\[
D^jf_e(S_\lambda(X,Z))[S_\lambda h,\ldots,S_\lambda h]
=
D^jf_e(X,Z)[h,\ldots,h],
\quad j=1,2,3.
\]
The bilinear form $G$ and the trilinear form $\mathcal T_f$ in Lemma~\ref{lem:projective} transform accordingly under $S_\lambda$, and the inequality \eqref{eq:projectiveSC} is invariant under this change of coordinates.
Given any $(X,Z)\in\D_e$, choose $\lambda=X^{-1}$. Then
$
S_{X^{-1}}(X,Z)=\left(1,\frac{Z}{X^c}\right).
$
Writing
$
t:=\frac{Z}{X^c},
$
the domain condition $|Z|<X^c$ is exactly $|t|<1$. Hence it suffices to verify \eqref{eq:projectiveSC} at points
\begin{equation}\label{eq:normalizedpoint}
(X,Z)=(1,t),
\quad
|t|<1.
\end{equation}
The symmetry $f_e(X,-Z)=f_e(X,Z)$ allows us to assume $0\le t<1$.

For a direction $(h_X,h_Z)$ at $(1,t)$, we introduce coordinates $(p,r)$ by
\begin{equation}\label{eq:dirchange}
h_X=p,
\quad
h_Z=ctp+r.
\end{equation}

Under the scaling $(X,Z)\mapsto(\lambda X,\lambda^c Z)$, the point $(1,t)$ is mapped to $(\lambda,\lambda^ct)$, whose derivative at $\lambda=1$ is $(1,ct)$. Hence every direction can be written uniquely as $(h_X,h_Z)=p(1,ct)+r(0,1)$. This choice also leads to a convenient polynomial form of the projective inequality in Lemma~\ref{lem:projective}.
\subsection{Quadratic and cubic forms}

Direct differentiation of \eqref{eq:fdef}, followed by the substitutions \eqref{eq:normalizedpoint}--\eqref{eq:dirchange}, yields
\begin{equation}\label{eq:GQ}
G[h,h]
=
\frac{9-e^2}{36A^2(1-t^2)^2}\,Q(p,r),
\end{equation}
where $A:=e^2+3$, $h=(h_X,h_Z)$, and $Q(p,r)$ is the quadratic form
\begin{equation}\label{eq:Q}
\begin{aligned}
Q(p,r)
={}&
(1-t^2)A(A-4e^2t^2)p^2
+8etA(1-t^2)pr 
+4(4e^2t^2+3e^2+9)r^2.
\end{aligned}
\end{equation}
The corresponding cubic form is
\begin{equation}\label{eq:CP}
\mathcal T_f[h,h,h]
=
-\frac{e(9-e^2)}{54A^3(1-t^2)^3}\,P(p,r),
\end{equation}
where the cubic polynomial is
\begin{equation}\label{eq:P}
\begin{aligned}
P(p,r)
={}&
A^3(1-t^2)^3p^3 
+12etA(1-t^2)
\bigl[A+2t^2(e^2-3)\bigr]p^2r \\
&
-12A(1-t^2)
\bigl[3A+t^2(5e^2-9)\bigr]pr^2 
-16et
\bigl[9A+t^2(5e^2-9)\bigr]r^3.
\end{aligned}
\end{equation}

Lemma~\ref{lem:projective} requires $G\succ0$. By \eqref{eq:GQ}, this is equivalent to the positive definiteness of the quadratic form $Q$.

\begin{lemma}\label{lem:Qpd}
For $0<e<1$, $0\le t<1$, and $(p,r)\ne(0,0)$,
\[
Q(p,r)>0.
\]
\end{lemma}

\begin{proof}
Set $u:=e^2$ and $s:=t^2$.
Thus $0<u<1$ and $0\le s<1$. Let $M$ be the symmetric matrix satisfying
\[
Q(p,r)=\begin{pmatrix}p&r\end{pmatrix}M\begin{pmatrix}p\\r\end{pmatrix}.
\]
Its lower-right entry is positive, and its determinant factors as
\begin{equation}\label{eq:detM}
\det M
=
12(1-s)(u+3)R(u,s),
\end{equation}
where $
R(u,s)
=
(u+3)^2-4u(u+3)s+4u(1-u)s^2.
$
Furthermore,
\[
\frac{\partial R}{\partial s}
=
-4u(u+3)+8u(1-u)s
\le
-4u(3u+1)<0.
\]
Hence $R(u,s)>R(u,1)$ for $s<1$ and
$
R(u,1)=(1-u)(7u+9)>0.
$
Thus $\det M>0$. Together with $M_{22}>0$, this implies that  $M\succ 0$.
\end{proof}

Using \eqref{eq:nu-e}, \eqref{eq:GQ}, and \eqref{eq:CP}, inequality \eqref{eq:projectiveSC} in Lemma~\ref{lem:projective} is equivalent to
\begin{equation}\label{eq:Deltaineq}
\Delta(p,r):=Q(p,r)^3-P(p,r)^2\ge0.
\end{equation}
Indeed,
\begin{equation}\label{eq:masteridentity}
4\gamma^2G[h,h]^3-\mathcal T_f[h,h,h]^2
=
\frac{e^2(9-e^2)^2}{2916A^6(1-t^2)^6}\,\Delta(p,r).
\end{equation}
It remains to prove \eqref{eq:Deltaineq}.

\subsection{Nonnegativity of the binary sextic} 
We first treat the case $t=0$ directly. For $0<t<1$, we show that the sextic has fixed degree and nonvanishing discriminant throughout the parameter domain. Its number of real roots is therefore constant, and a single exact Sturm computation determines the root count on the whole domain.

\begin{lemma}\label{lem:t0}
If $t=0$, then
\begin{equation}\label{eq:t0factor}
\Delta(p,r)
=
108r^2A^3\bigl(Ap^2-4r^2\bigr)^2
\ge0.
\end{equation}
\end{lemma}

\begin{proof}
Substitute $t=0$ into \eqref{eq:Q}--\eqref{eq:P}. One obtains
\[
Q=A^2p^2+12Ar^2,
\quad
P=A^3p^3-36A^2pr^2.
\]
Expanding $Q^3-P^2$ and factoring the result gives \eqref{eq:t0factor}. Since $A=e^2+3>0$, every factor on the right-hand side is nonnegative.
\end{proof}

Since $\Delta$ is homogeneous of degree six, it follows that
$
\Delta(p,r)=p^6\Delta(1,r/p)
$
whenever $p\ne0$. For $0<t<1$, define
\begin{equation}\label{eq:ddef}
d_{e,t}(\xi):=\Delta(1,\xi).
\end{equation}
Since $Q(1,\xi)$ is quadratic and $P(1,\xi)$ is cubic, $d_{e,t}$ has degree at most six. Its leading coefficient is determined by the highest-degree terms of $Q(1,\xi)$ and $P(1,\xi)$. From \eqref{eq:Q},
\[
[\xi^2]Q(1,\xi)
=
4(4e^2t^2+3e^2+9),
\]
and from \eqref{eq:P},
\[
[\xi^3]P(1,\xi)
=
-16et\bigl[9A+t^2(5e^2-9)\bigr].
\]
Therefore
\begin{align}
[\xi^6]d_{e,t}
&=
\bigl([\xi^2]Q(1,\xi)\bigr)^3
-
\bigl([\xi^3]P(1,\xi)\bigr)^2\notag\\
&=
64\Bigl[
(4e^2t^2+3e^2+9)^3
-4e^2t^2\bigl(9(e^2+3)+t^2(5e^2-9)\bigr)^2
\Bigr].\label{eq:leading-before-L}
\end{align}
Now set $u=e^2$ and $s=t^2$ and define the auxiliary polynomial
\begin{equation}\label{eq:Lexplicit}
L(u,s)
:=
\frac19\Bigl[
(4us+3u+9)^3
-4us\bigl(9(u+3)+s(5u-9)\bigr)^2
\Bigr].
\end{equation}
Consequently, the leading coefficient is
\begin{equation}\label{eq:leading}
[\xi^6]d_{e,t}=576L(u,s).
\end{equation}
To prove that $d_{e,t}$ has degree six, it remains to show that $L(u,s)>0$ for $0<u,s<1$.

The expansion of $L$ in the power basis is
\begin{align}
L(u,s)
={}&81+81u+27u^2+3u^3
-216us-144u^2s-24u^3s\notag\\
&+216us^2-24u^3s^2
-36us^3+40u^2s^3-4u^3s^3.\label{eq:Lpower}
\end{align}
Because the coefficients in the power basis have mixed signs, we rewrite $L$ in the Bernstein basis introduced in Subsection~\ref{subsec:poly-prelim}. Applying \eqref{eq:power-to-bernstein} in both variables with degree $3$ yields
\begin{equation}\label{eq:Lbernstein}
L(u,s)
=
\sum_{i=0}^3\sum_{j=0}^3
\ell_{ij}B_i^{(3)}(u)B_j^{(3)}(s),
\end{equation}
where the coefficient matrix is
\begin{equation}\label{eq:Lmatrix}
(\ell_{ij})_{i,j=0}^3
=
\begin{pmatrix}
81&81&81&81\\
108&84&84&96\\
144&80&64&256/3\\
192&64&0&0
\end{pmatrix}.
\end{equation}
All entries are nonnegative. In fact, using only the first row and the partition-of-unity identity,
\[
L(u,s)
\ge
81B_0^{(3)}(u)\sum_{j=0}^3B_j^{(3)}(s)
=
81(1-u)^3.
\]
Hence
\begin{equation}\label{eq:Lpositive}
L(u,s)>0
\quad
\text{for }0<u,\,s<1.
\end{equation}
In particular, the leading coefficient in \eqref{eq:leading} never vanishes,
so $d_{e,t}$ has degree exactly six throughout the open parameter square,
as required in Lemma~\ref{lem:root-continuation}.

\subsubsection{Discriminant factorization}

We now show that $d_{e,t}$ has no multiple complex roots. Since the leading coefficient is nonzero, Lemma~\ref{lem:disc-background} shows that this is
equivalent to
$
\Disc_\xi d_{e,t}\ne0.
$
Since $d_{e,t}$ has degree six, the definition of the discriminant gives
\[
\Disc_\xi d_{e,t}
=
-\frac{\Res_\xi(d_{e,t},\partial_\xi d_{e,t})}{[\xi^6]d_{e,t}},
\]
Computing the corresponding Sylvester determinant and factoring the resulting polynomial yields

\begin{equation}\label{eq:discd}
\begin{aligned}
\Disc_{\xi}d_{e,t}
={}&
2^{44}3^{27}s^3(u-1)^3(u+3)^{15}(s-1)^{15} 
\Pi_1(u,s)^3\Pi_2(u,s),
\end{aligned}
\end{equation}
where the two polynomials $\Pi_1$ and $\Pi_2$ are listed explicitly in Appendix~\ref{app:polynomials}.

On the open square $0<u,s<1$, the factors
$s^3$, $(u-1)^3$, $(u+3)^{15}$, and $(s-1)^{15}$
are all nonzero. Consequently, the discriminant can vanish only if either $\Pi_1$ or $\Pi_2$ vanishes. It therefore suffices to determine the signs of $\Pi_1$ and $\Pi_2$ throughout $(0,1)^2$.

\subsubsection{Sign of \texorpdfstring{$\Pi_1$}{Pi1}}
We determine the sign of $\Pi_1$ using the trace-free part of the symmetric third-order tensor associated with $P$. With repeated indices understood to be summed, write
\[
Q(h)=M_{ij}h_i h_j,
\quad
P(h)=T_{ijk}h_i h_j h_k,
\]
where $M\succ0$ and $T$ is symmetric.  We use the standard correspondence between homogeneous polynomials and symmetric tensors~\cite{ComonGolubLimMourrain2008,Sharafutdinov1994}.

To see why divisibility by $Q$ is relevant, suppose that $P(h)=Q(h)\ell(h)$ for some linear form $\ell$. Then
\begin{equation}\label{eq:pure-trace-degeneration}
\Delta(h)
=
Q(h)^3-P(h)^2
=
Q(h)^2\bigl(Q(h)-\ell(h)^2\bigr).
\end{equation}
Since $M\succ0$, the polynomial $Q(1,\xi)$ has two distinct complex roots. If $P=Q\ell$, then \eqref{eq:pure-trace-degeneration} shows that these are multiple roots of $d_{e,t}(\xi)=\Delta(1,\xi)$. Hence $\Disc_\xi d_{e,t}=0$ by Lemma~\ref{lem:disc-background}.

Motivated by this observation, write $M^{-1}=(M^{ij})$ and define the trace of $T$ with respect to $M$ by
\begin{equation}\label{eq:tensor-trace}
\tau_i:=M^{jk}T_{ijk}.
\end{equation}
Its action on $h$ is then given by $\tau(h):=\tau_i h_i$. Following the standard trace-free decomposition of symmetric third-order tensors \cite[Section~3]{TothTuryshev2022}, define
\begin{equation}\label{eq:Htensor}
H_{ijk}
:=
T_{ijk}
-\frac14
\bigl(
M_{ij}\tau_k+M_{ik}\tau_j+M_{jk}\tau_i
\bigr).
\end{equation}
The coefficient $1/4$ is chosen so that $H$ is
trace-free with respect to $M$. Indeed,
\[
\begin{aligned}
M^{jk}H_{ijk}
&=
\tau_i
-\frac14
\bigl(
M^{jk}M_{ij}\tau_k
+M^{jk}M_{ik}\tau_j
+M^{jk}M_{jk}\tau_i
\bigr) \\
&=
\tau_i-\frac14(\tau_i+\tau_i+2\tau_i)
=0.
\end{aligned}
\]

Let $H(h):=H_{ijk}h_i h_j h_k$ be the corresponding cubic form.
Contracting \eqref{eq:Htensor} with $h_i h_j h_k$ gives
\begin{equation}\label{eq:Hpoly}
H(h)
=
P(h)-\frac34Q(h)\tau(h).
\end{equation}
In particular, if $H=0$, then
$
P(h)=\frac34Q(h)\tau(h),
$
so $P$ is divisible by $Q$. In this case,
\begin{equation}\label{eq:H-zero-degeneration}
\Delta(h)
=
Q(h)^2
\left(
Q(h)-\frac9{16}\tau(h)^2
\right),
\end{equation}
and the sextic discriminant vanishes for the reason described above.

It remains to connect $H$ with the factor $\Pi_1$. Since $M\succ0$, its inverse defines a positive-definite inner product on symmetric third-order tensors. We use
\begin{equation}\label{eq:tensor-M-norm}
\langle S,T\rangle_M
:=
M^{ia}M^{jb}M^{kc}S_{ijk}T_{abc},
\quad
\|T\|_M^2:=\langle T,T\rangle_M.
\end{equation}
This is positive definite because, in an $M$-orthonormal basis, it is the usual sum of products of tensor entries.

Substituting the coefficients of $Q$ and $P$ into
\eqref{eq:Htensor} and simplifying the contraction in \eqref{eq:tensor-M-norm} gives
\begin{equation}\label{eq:Hnormidentity}
\|H\|_M^2
=
M^{ia}M^{jb}M^{kc}H_{ijk}H_{abc}
=
-\frac{\Pi_1(u,s)}{R(u,s)^3}.
\end{equation}
Using \eqref{eq:tensor-trace} and substituting the coefficients of $Q$ and $P$ gives
$
\tau_1=M^{11}T_{111}+2M^{12}T_{112}+M^{22}T_{122}=0.
$
It then follows from \eqref{eq:Htensor} that
\begin{equation}\label{eq:H111}
H_{111}=T_{111}-\frac{3}{4}M_{11}\tau_1=A^3(1-s)^3>0
\quad
\text{for }0<s<1.
\end{equation}
Thus, $H\neq0$, and the positive definiteness in
\eqref{eq:tensor-M-norm} yields $\|H\|_M^2>0$.
The proof of Lemma~\ref{lem:Qpd} also gives $R(u,s)>0$ for
$0<u,s<1$. Combining these facts with \eqref{eq:Hnormidentity} implies
\begin{equation}\label{eq:pi1negative}
\Pi_1(u,s)<0
\quad
\text{for }0<u,\,s<1.
\end{equation}

\subsubsection{Sign of \texorpdfstring{$\Pi_2$}{Pi2}}\label{signofPi}  
For $\Pi_2$, we use a continuation argument instead.

\begin{lemma}\label{lem:compactcontinuation}
Let $I\subset\R$ be an interval and let $p(u,s)$ be a real polynomial. Suppose, for every $u\in I$,
\[
p(u,0)>0,
\quad
p(u,1)>0,
\quad
\Disc_s p(u,\cdot)\ne0.
\]
If $p(u_0,s)>0$ for all $s\in[0,1]$ at one point $u_0\in I$, then
\[
p(u,s)>0
\quad
\text{for all }(u,s)\in I\times[0,1].
\]
\end{lemma}

\begin{proof}
Fix any $u_1\in I$. It is enough to prove positivity on the compact line segment joining $u_0$ and $u_1$. Parametrize this segment by $u(\theta)=(1-\theta)u_0+\theta u_1$, $0\le\theta\le1$, and define
\[
m(\theta):=\min_{0\le s\le1}p(u(\theta),s).
\]
Since $p$ is continuous and the interval in $s$ is compact, the minimum is attained. Moreover, the minimum of a continuous function over a fixed compact interval depends continuously on the remaining parameters, so $m$ is continuous.
We have $m(0)>0$. If $m(1)\le0$, let
\[
\theta_*:=\inf\{\theta\in[0,1]:m(\theta)\le0\}.
\]
Continuity gives $m(\theta_*)=0$. Choose $s_*\in[0,1]$ attaining this minimum. The endpoint hypotheses $p(u,0)>0$ and $p(u,1)>0$ force $s_*\in(0,1)$. Since $s_*$ is an interior minimizer,
\[
p(u(\theta_*),s_*)=0,
\quad
\partial_s p(u(\theta_*),s_*)=0.
\]
Thus $p(u(\theta_*),\cdot)$ and its derivative have a common root, so Lemma~\ref{lem:disc-background} gives
\[
\Disc_s p(u(\theta_*),\cdot)=0,
\]
contrary to the assumption. Hence $m(1)>0$. Since $u_1$ was arbitrary, the conclusion follows.
\end{proof}

For $0<u<1$, the boundary values of $\Pi_2$ are
\begin{equation}\label{eq:pi2boundary}
\Pi_2(u,0)=4(u+3)^9>0,
\quad
\Pi_2(u,1)=u(1-u)^6(7u+9)^2>0.
\end{equation}
Its discriminant with respect to $s$ factors as
\begin{equation}\label{eq:discpi2}
\begin{aligned}
\Disc_s\Pi_2
={}&
-2^{12}u^8(u-1)^{24}(u+3)^{72}(u+7)^6 \\
&\times(2u-1)^3(7u+9)^3A_6(u)^3A_7(u),
\end{aligned}
\end{equation}
where the auxiliary polynomials $A_6$ and $A_7$ are given in Appendix~\ref{Positivity of A}. The same appendix gives exact Bernstein coefficients certifying
\[
A_6(u)>0,
\quad
A_7(u)>0
\quad \text{for }
0\le u\le1.
\]
It follows from \eqref{eq:discpi2} that all factors are nonzero for $0<u<1$ except $(2u-1)^3$, which vanishes at $u=1/2$. Therefore
\begin{equation}\label{eq:discpi2nonzero}
\Disc_s\Pi_2(u,\cdot)\ne0
\end{equation}
for the two connected intervals
$
0<u<\frac12
$ and
$
\frac12<u<1.
$
Lemma~\ref{lem:compactcontinuation} can be applied separately on these intervals. It does not apply on an interval containing $u=1/2$, since $\Disc_s\Pi_2(u,\cdot)$ vanishes at $u=1/2$.

It remains to verify positivity at one base point in each open interval and at the exceptional value $u=1/2$. The exact algebraic decompositions in Appendix~\ref{Base-point checks} show that
\[
\Pi_2\left(\frac14,s\right)>0,
\quad
\Pi_2\left(\frac34,s\right)>0,
\quad
\Pi_2\left(\frac12,s\right)>0
\]
for every $0\le s\le1$.
Applying Lemma~\ref{lem:compactcontinuation} with $u_0=1/4$ gives
\begin{equation}\label{eq:pi2firsthalf}
\Pi_2(u,s)>0
\quad
\text{for }0<u<\frac12,\ 0\le s\le1,
\end{equation}
and applying it with $u_0=3/4$ gives
\begin{equation}\label{eq:pi2secondhalf}
\Pi_2(u,s)>0
\quad
\text{for }\frac12<u<1,\ 0\le s\le1.
\end{equation}
Together with the verification at $u=1/2$, these estimates yield
\begin{equation}\label{eq:pi2positive}
\Pi_2(u,s)>0
\quad
\text{for }0<u<1,\ 0\le s\le1.
\end{equation}

\subsubsection{Strict positivity of \texorpdfstring{$d_{e,t}$}{d(e,t)}}
By \eqref{eq:discd}, \eqref{eq:pi1negative}, and \eqref{eq:pi2positive}, we have
\begin{equation}\label{eq:discdnonzero}
\Disc_{\xi}d_{e,t}\ne0
\quad
\text{for }0<e,\,t<1.
\end{equation}
Together with \eqref{eq:Lpositive}, the hypotheses of Lemma~\ref{lem:root-continuation} are satisfied on the connected parameter square $(0,1)^2$.  The number of real roots of $d_{e,t}$ is therefore constant there. We determine it at $(e,t)=(\frac13,\frac12)$. 
Since multiplication by a positive constant does not affect the real roots, define
\begin{equation}\label{eq:basepoly}
\begin{aligned}
d_0(\xi)
:=
\frac{729}{2}d_{1/3,\,1/2}(\xi)
={}&
17914752\xi^6-5322240\xi^5-3190320\xi^4
+3596208\xi^3 \\
&+3949008\xi^2+94668\xi+66199.
\end{aligned}
\end{equation}
A Sturm chain for $d_0$, with its members individually rescaled by positive constants, is listed in Appendix~\ref{app:sturm}. 
The signs of the leading terms at $-\infty$ and $+\infty$ are, respectively,
$(+,-,+,+,+,-,-)$ and $(+,+,+,-,+,+,-)$. Since both sequences have three sign variations,
Lemma~\ref{lem:sturm-background} gives
\[
V(-\infty)-V(+\infty)=3-3=0,
\]
so $d_0$ has no real root.

Because $d_0$ is continuous and never vanishes on $\R$, its sign is constant. Since $d_0(0)=66199>0$, it follows that
\[
d_0(\xi)>0
\quad
\text{for every }\xi\in\R.
\]
Lemma~\ref{lem:root-continuation} now shows that every $d_{e,t}$ has no real roots. Since $d_{e,t}$ has degree six with positive leading coefficient, it tends to $+\infty$ as $|\xi|\to\infty$. A continuous real polynomial with no real root cannot change sign; hence it is everywhere positive. Therefore
\begin{equation}\label{eq:dpositive}
d_{e,t}(\xi)>0
\quad
\text{for }0<e,t<1,\ \xi\in\R.
\end{equation}

The polynomial $\Delta(p,r)=Q(p,r)^3-P(p,r)^2$ is homogeneous of degree six because $Q$ is homogeneous of degree two and $P$ of degree three. Hence, if $p\ne0$,
\[
\Delta(p,r)
=p^6\Delta\left(1,\frac rp\right)
=p^6d_{e,t}(r/p)>0.
\]
If $p=0$ and $r\ne0$, the coefficient $[\xi^6]d_{e,t}$ is exactly the coefficient of $r^6$ in the binary sextic $\Delta(p,r)$, so \eqref{eq:leading} gives
\[
\Delta(0,r)=576L(u,s)r^6>0.
\]
Combining this with Lemma~\ref{lem:t0}, we have proved
\begin{equation}\label{eq:Deltaglobal}
\Delta(p,r)\ge0
\end{equation}
for every $0<e<1$, $0\le t<1$, and every $(p,r)$, with strict inequality for nonzero directions whenever $t>0$.

These properties allow us to invoke Lemma~\ref{lem:projective} and obtain the desired result.
\begin{theorem}\label{thm:upperbound}
For $0<c\le1/2$, the function $F_c$ in \eqref{eq:Fdef} is a $3/(1+c)$-LHSCB on $\Kpow_c$.
\end{theorem}
\begin{proof}
Subsection~\ref{subsec:4.1} establishes logarithmic homogeneity and the barrier property of $F_c$.
By \eqref{eq:Deltaglobal} and \eqref{eq:masteridentity},
inequality~\eqref{eq:projectiveSC} holds. Lemma~\ref{lem:Qpd} gives $G\succ0$. Lemma~\ref{lem:projective} therefore implies that $F_c$ has positive-definite Hessian and satisfies the standard self-concordance inequality for every $0<c<1/2$. For $c=1/2$, the barrier $F_c$ takes the form
\[
F_{1/2}(x,y,z)
=-\log x-\log y-\log\left(1-\frac{z^2}{xy}\right)
=-\log(xy-z^2),
\]
which is the standard LHSCB with parameter $2$ for the
rotated second-order cone. Since $3/(1+1/2)=2$, the proof is complete.
\end{proof}
\section{The optimal barrier parameter}\label{sec5}
We can now combine the lower and upper bounds.

\begin{proof}[Proof of Theorem~\ref{thm:main}]
For $0<c\le1/2$, Theorem~\ref{thm:lowerbound} gives $\nu_{\rm opt}(\Kpow_c)\ge3/(1+c)$, while Theorem~\ref{thm:upperbound} shows that the explicit barrier $F_c$ in \eqref{eq:Fdef} attains this value. Hence equality holds.
\end{proof}

\begin{remark}[Sharpness of the lower-bound argument]\label{rem:lower-sharp}
For $0<c<1/2$, consider the explicit barrier $F_c$ at the optimal
parameter $\nu_c=3/(1+c)$. Using the notation from the proof of
Theorem~\ref{thm:lowerbound}, the minimizer in the $z$-variable on the
section $x+y=1$ is $\zeta(r)=0$. Set
$
\varrho_c:=\nu_c-1=\frac{2-c}{1+c}.
$
For this barrier,
\[
\widehat q(r)
=
\frac{(1+c)(2-c)}
     {9r^2(1-r)^2},
\quad
k(r)
=
\frac{(1+c)(2-c)}
     {3(c^2-c+1)r^{2c}(1-r)^{2(1-c)}}.
\]
Consequently,
$
\kappa(r)
=
k(r)w(r)^2
=
\frac{(1+c)(2-c)}
     {3(c^2-c+1)}
$
is constant. Moreover,
\[
\ell_-=\sqrt{\varrho_c},
\quad
\ell_+=\frac{1}{\sqrt{\varrho_c}},
\quad
\beta=
\frac{2(\varrho_c-1)}{\sqrt{\varrho_c}}.
\]
Thus the upper bound for $\beta$, the upper bound for $\ell_-$, and the
lower bound for $\ell_+$ used in \eqref{eq:drift-bound} are all attained
with equality. In particular,
\[
\frac{\text{d}}{\text{d}\sigma}\log\kappa
=
0
=
\frac{2}{\sqrt{\varrho_c}}
\bigl((1+c)\varrho_c-(2-c)\bigr).
\]
Hence the estimate \eqref{eq:drift-bound} is sharp at the optimal parameter $\nu_c$.
\end{remark}

The symmetry in Lemma~\ref{lem:exponent-symmetry} extends the result to the remaining exponents.
\begin{corollary}[Full exponent range]\label{cor:full-range}
For every $0<c<1$,
\[
\nu_{\rm opt}(\Kpow_c)
=
\frac{3}{1+\min\{c,1-c\}}.
\]
If $c\le1/2$, the optimum is attained by $F_c(x,y,z)$ in \eqref{eq:Fdef}. If $c>1/2$, it is attained by $F_{1-c}(y,x,z)$.
\end{corollary}

\begin{proof}
The case $c\le1/2$ is Theorem~\ref{thm:main}. If $c>1/2$, Lemma~\ref{lem:exponent-symmetry} maps $\Kpow_c$ linearly onto $\Kpow_{1-c}$ by exchanging $x$ and $y$, and $0<1-c<1/2$. Applying Theorem~\ref{thm:main} to exponent $1-c$ proves both the parameter formula and the stated barrier.
\end{proof}

\section{Conclusion}\label{sec6}

We have determined the optimal barrier parameter of the
three-dimensional power cone. For every $0<c<1$, the optimal  parameter is
$
\nu_{\mathrm{opt}}(\Kpow_c)
=
\frac{3}{1+\min\{c,1-c\}},
$
and the optimum is attained by an explicit barrier.  At $c=1/2$, the formula recovers the optimal parameter $2$ of the second-order cone. The lower-bound argument applies to arbitrary LHSCBs, and the explicit optimal barrier attains equality in its key estimates, confirming the sharpness of the bound. The matching upper bound is obtained by reducing the self-concordance inequality to the nonnegativity of a binary sextic and proving this property by exact algebraic arguments. A natural direction for further work is to determine whether the optimal barrier parameters can be obtained for higher-dimensional generalized power cones.

\appendix

\section{Exact polynomial factors}\label{app:polynomials}
The two factors in \eqref{eq:discd} have the expansions
\[
\Pi_1(u,s)=\sum_{j=0}^6a_j(u)s^j,
\quad
\Pi_2(u,s)=\sum_{j=0}^9b_j(u)s^j.
\]
Their coefficients are given in the following compact form:
\begingroup
\setlength{\jot}{2pt}
\allowdisplaybreaks
\begin{align*}
a_0(u)&=-4(u+3)^6,\\
a_1(u)&=12(u+3)^6,\\
a_2(u)&=12(u+3)^4(3u^2-10u-9),\\
a_3(u)&=-(u+3)^3(23u^3+294u^2-465u-108),\\
a_4(u)&=12u(u-1)(u+3)^2(2u^2-u+63),\\
a_5(u)&=-12u(u-3)(u-1)^2(u+3)(u+9),\\
a_6(u)&=16u^2(u-1)^2(u^2+2u+9);\\[2pt]
b_0(u)&=4(u+3)^9,\\
b_1(u)&=-(u+3)^8(u+72),\\
b_2(u)&=(u+3)^7(u^2-45u+576),\\
b_3(u)&=(u+3)^6(13u^3-42u^2+777u-2700),\\
b_4(u)&=-4(u+3)^5(5u^4+27u^3-177u^2+1269u-2052),\\
b_5(u)&=4(u+3)^4(6u^5+21u^4+117u^3-1237u^2+4473u-4212),\\
b_6(u)&=4(u+3)^3(6u^6-81u^5+111u^4-662u^3+4392u^2-9369u+5859),\\
b_7(u)&=4(u-1)(u+3)^2
\bigl(u^6-50u^5+340u^4-492u^3+2049u^2-6426u+5346\bigr),\\
b_8(u)&=-24(u-1)^2(u+3)
\bigl(u^5-22u^4+82u^3-108u^2+405u-486\bigr),\\
b_9(u)&=4(u-9)(u-1)^3(9u^3-19u^2+27u-81).
\end{align*}
\endgroup

\section{Auxiliary polynomial calculations}\label{app:bernstein}
This appendix contains the auxiliary polynomial calculations used in Section~\ref{sec4}. 

\subsection{Positivity of \texorpdfstring{$A_6$ and $A_7$}{A6 and A7}}\label{Positivity of A}

The auxiliary factors in \eqref{eq:discpi2} are
\begin{align*}
A_6(u)
={}&1612u^6-7822u^5+35755u^4-62928u^3
+62694u^2-18954u+2187,\\
A_7(u)
={}&46301u^7-430305u^6+715218u^5+2010910u^4 \\
&-7240815u^3+31031667u^2-12772080u+1679616.
\end{align*}
To certify positivity on $[0,1]$, we divide the interval into $[0,1/2]$ and $[1/2,1]$. The affine substitutions $u=v/2$ and $u=(1+v)/2$ map $[0,1]$ onto these two subintervals, respectively. For $A_6(v/2)$, the Bernstein coefficients are
\[
\left(
2187,
\frac{1215}{2},
\frac{729}{10},
\frac{1899}{10},
\frac{171403}{240},
\frac{48309}{32},
\frac{40527}{16}
\right).
\]
For $A_6((1+v)/2)$, they are
\[
\left(
\frac{40527}{16},
\frac{113799}{32},
\frac{1153753}{240},
\frac{31214}{5},
7912,
\frac{29824}{3},
12544
\right).
\]
Hence $A_6(u)>0$ on $[0,1]$.

For $A_7(v/2)$, the Bernstein coefficients are
\[
\left(
1679616,
\frac{5371272}{7},
\frac{897831}{4},
\frac{1408707}{56},
\frac{1177189}{8},
\frac{64156205}{112},
\frac{577617375}{448},
\frac{292871875}{128}
\right).
\]
For $A_7((1+v)/2)$, they are
{\footnotesize
\[
\left(
\frac{292871875}{128},
\frac{736242875}{224},
\frac{1023180785}{224},
\frac{42856724}{7},
\frac{55624276}{7},
\frac{70273152}{7},
\frac{86818176}{7},
15040512
\right).
\]
}
Hence $A_7(u)>0$ on $[0,1]$.

\subsection{Base-point checks for \texorpdfstring{$\Pi_2$}{Pi2}}\label{Base-point checks}

We give the exact decompositions used for the three base-point checks in Subsection~\ref{signofPi}. In each case, set $q=1-s$.

For $u=1/4$,
\begin{equation}\label{eq:pi2quarter}
\begin{aligned}
2^{18}\Pi_2\left(\frac14,1-q\right)
={}&27h_1(q)
+6450327560q^4+12966173284q^5 \\
&+12808009492q^6+6652881036q^7 \\
&+1866026592q^8+291453120q^9,
\end{aligned}
\end{equation}
where
\begin{equation}\label{eq:h1}
h_1(q)
=
50379285q^3+1563956q^2-766260q+49923.
\end{equation}
Its discriminant is
\begin{equation}\label{eq:disch1}
\Disc h_1
=-133708464211246463833953147<0.
\end{equation}
By Lemma~\ref{lem:disc-background}, $h_1$ has exactly one real root. Its leading coefficient is positive and $h_1(0)=49923>0$, so this root is negative. Hence $h_1(q)>0$ for $q\ge0$. All remaining terms on the right-hand side of \eqref{eq:pi2quarter} are nonnegative for $q\ge0$. Therefore
\[
\Pi_2\left(\frac14,s\right)>0
\quad
\text{for }0\le s\le1.
\]

For $u=3/4$,
\begin{equation}\label{eq:pi2threequarter}
\begin{aligned}
2^{18}\Pi_2\left(\frac34,1-q\right)
={}&27h_2(q)
+54014298552q^4+69669314028q^5 \\
&+25199957628q^6+3569396868q^7 \\
&+234936288q^8+9142848q^9,
\end{aligned}
\end{equation}
where
\begin{equation}\label{eq:h2}
h_2(q)
=
40628031q^3-747684q^2-14364q+361.
\end{equation}
Its discriminant is
\begin{equation}\label{eq:disch2}
\Disc h_2
=-1772210425734147611691<0.
\end{equation}
Again Lemma~\ref{lem:disc-background} gives exactly one real root. Since the leading coefficient is positive and $h_2(0)=361>0$, that root is negative. Hence $h_2(q)>0$ for $q\ge0$, and all higher-degree terms in \eqref{eq:pi2threequarter} are nonnegative there. Therefore
\[
\Pi_2\left(\frac34,s\right)>0
\quad
\text{for }0\le s\le1.
\]

Finally, at $u=1/2$,
\begin{equation}\label{eq:pi2half}
\Pi_2\left(\frac12,s\right)
=-\frac{(2s-7)^2}{512}R_7(s),
\end{equation}
where
\begin{equation}\label{eq:R7}
\begin{aligned}
R_7(s)
={}&
38692s^7-539840s^6+3911572s^5-14174132s^4 \\
&+27167315s^3-28286181s^2+15176721s-3294172.
\end{aligned}
\end{equation}
With $q=1-s$,
\begin{equation}\label{eq:R7q}
-R_7(1-q)
=h_0(q)+1359652q^4+1485064q^5+268996q^6+38692q^7,
\end{equation}
where
\begin{equation}\label{eq:h0}
h_0(q)=143927q^3-1624q^2-560q+25.
\end{equation}
Its discriminant is
\begin{equation}\label{eq:disch0}
\Disc h_0=-188304545926875<0.
\end{equation}
By Lemma~\ref{lem:disc-background}, $h_0$ has exactly one real root. Since its leading coefficient and $h_0(0)=25$ are positive, this root is negative. Hence $h_0(q)>0$ for $q\ge0$. All remaining terms in \eqref{eq:R7q} have positive coefficients, so $-R_7(1-q)>0$ for $q\ge0$. Therefore $R_7(s)<0$ for $0\le s\le1$. Since $(2s-7)^2>0$ on $[0,1]$, \eqref{eq:pi2half} gives
\[
\Pi_2\left(\frac12,s\right)>0
\quad
\text{for }0\le s\le1.
\]

\section{A Sturm chain for \texorpdfstring{$d_0$}{d0}}\label{app:sturm}
For the polynomial $d_0$ in \eqref{eq:basepoly}, start from $S_0=d_0$ and $S_1=d_0'$ and repeatedly take the negative Euclidean remainder, as in Lemma~\ref{lem:sturm-background}. The following list is the resulting Sturm chain after multiplication of individual members by positive rational constants. Those rescalings keep every sign unchanged and therefore leave the Sturm count invariant:
\allowdisplaybreaks
\begin{align*}
S_0(\xi)={}&
17914752\xi^6-5322240\xi^5-3190320\xi^4+3596208\xi^3
+3949008\xi^2\\&+94668\xi+66199,\\
S_1(\xi)={}&
8957376\xi^5-2217600\xi^4-1063440\xi^3+899052\xi^2+658168\xi+7889,\\
S_2(\xi)={}&
2850379920\xi^4-3760659672\xi^3-6046460336\xi^2
-320056730\xi-148801387,\\
S_3(\xi)={}&
-330462921868609056\xi^3
-240471376030379788\xi^2\\
&
-23796321104735200\xi
-5496922181504621,\\
S_4(\xi)={}&
138655679724599952107248724\xi^2
-3642291457427814879633064\xi\\
&
+3576940775423415451698955,\\
S_5(\xi)={}&
3292010764698075971634963590183272\xi\\
&
-140414792326549871129269285851577,\\
S_6(\xi)={}&-1.
\end{align*}

\section*{Data availability statement}
The Python code used to verify the algebraic identities and symbolic computations appearing in this paper is available at \url{https://github.com/wenhfu/powercone-verification}.

\section*{Acknowledgments}
 Generative AI tools were used in checking some lengthy algebraic manipulations and in assisting with the preparation of computational verification scripts. All mathematical arguments, identities, and conclusions reported in this paper were independently checked by the authors. The authors assume responsibility for all content.

\bibliographystyle{siamplain}
\bibliography{references}
\end{document}